\documentclass[journal]{IEEEtran}
\usepackage[cmex10]{amsmath}
\usepackage{amsfonts}
\usepackage{amssymb}
\usepackage{graphicx}
\usepackage{verbatim}
\usepackage{array}
\usepackage{multirow}
\usepackage{dcolumn}
\usepackage{color}
\usepackage[noadjust]{cite}
\usepackage{url}
\usepackage{balance}
\usepackage{booktabs}
\usepackage{subfigure}
\usepackage{wasysym}
\usepackage[mathscr]{euscript}
\DeclareGraphicsExtensions{.eps}
\usepackage{setspace}
\usepackage{graphicx,epstopdf}
\usepackage{comment}
\usepackage{xcolor}
\usepackage{flushend}

\begin{document}
\title{Assessing Costs of Community Solar Integration via Optimal Distribution Grid Expansion}

\author{Miguel Heleno, \IEEEmembership{Senior Member,~IEEE,} Alan Valenzuela, Alexandre Moreira \IEEEmembership{Member,~IEEE}, Juan Pablo Carvallo.

\thanks{ M. Heleno, A. Valenzuela, A. Moreira and J. P. Carvallo are with the Lawrence Berkeley National Laboratory, Berkeley, CA, USA (e-mail: \mbox{\{MiguelHeleno, AlanValenzuela, AMoreira,JPCarvallo\}@lbl.gov}).}
}

\maketitle
\begin{abstract}
This paper presents a comprehensive analysis of distribution grid upgrades costs necessary to integrate Community Solar (CS) projects. The innovative methodology proposed for this quantification is based on incremental least-cost
expansion of the distribution system, encompassing both traditional and non-traditional grid upgrade strategies. Realistic infrastructure investment costs are obtained using a dataset of over 2,500 feeders, including various loading scenarios, and compared with empirical costs of integrating real CS projects. The results are then used to evaluate costs and deferrals from the consumers’ and project developers’ perspectives, assess the nature of the infrastructure upgrades and explore the potential benefits of a strategic siting of CS projects. This analysis is summarized in a set of regulatory and policy recommendations to support planning and valuation aspects related to CS.
\end{abstract}

\section*{Nomenclature}\label{Nomenclature}
    
The mathematical symbols used throughout this paper are classified below as follows.

\subsection*{Sets}

\begin{description} [\IEEEsetlabelwidth{50000000}\IEEEusemathlabelsep]
    \item[${\cal D}$] Set of indexes of representative days.
    
    \item[$H$] Set of indexes of storage units.
    
    \item[$H^{C}$] Set of indexes of candidate storage units.

    \item[${\cal L}$] Set of indexes of all line segments.

    \item[${\cal L}^F$] Set of indexes of fixed line segments.

    \item[${\cal L}^{FH}$] Set of indexes of feeder head.

    \item[${\cal L}^C$] Set of indexes of candidate line segments for capacity upgrade.

    \item[${\cal L}^{VR}$] Set of indexes of line segments with existing or candidate voltage regulators.

    \item[${\cal L}^{VR,E}$] Set of indexes of line segments with existing voltage regulators.

    \item[${\cal L}^{VR,C}$] Set of indexes of line segments with candidate voltage regulators.
    
    \item[$N$] Set of indexes of all buses.


    \item[$R^{line,options}_l$] Set of indexes of upgrade capacity options for line segment $l$.

    \item[$R^{Solar,C}$] Set of indexes of candidate solar units.

    \item[$R^{Solar,E}$] Set of indexes of existing solar units.

    \item[$R^{Subs}_n$] Set of indexes of substations at bus $n$. This set will have at most one index for each $n$.
    
    \item[$T$] Set of all time periods.
    
    
\end{description}

\vspace{-0.3cm}

\subsection*{Parameters}

\begin{description} [\IEEEsetlabelwidth{10000000}\IEEEusemathlabelsep]

    \item[$\beta^p_l,\beta^q_l$] Auxiliary parameters related to active and reactive power losses, respectively.


    \item[$\eta$] Round-trip efficiency of storage units.

    \item[$\varphi^{max}_l,\varphi^{min}_l$] Maximum and minimum turns ratio of voltage regulator.


    \item[$\psi$] Auxiliary parameter related to the reactive contribution capacity of storage units.

    \item[$\omega_w$] Weight of representative day $d$.

    \item[$C^{CS,crt}_t$] Cost of curtailing solar generation output.

    \item[$C^{I}$] Imbalance cost.

    \item[$C^{FH,inv}_l$] Equivalent annual investment cost to upgrade the feeder head transformer.

    \item[$C^{line,inv}_l$] Equivalent annual investment cost to upgrade line capacity.

    \item[$C^{VR,inv}_l$] Equivalent annual investment cost to updgrade line segment $l$ with an voltage regulator.
    
    \item[$C^{st,inv}_h$] Equivalent annual investment cost to install storage $h$.

    \item[${f}^{RTS}_{r,t,d}$] Capacity factor existing solar unit $r$.

    \item[${f}^{CS}_{r,t,d}$] Capacity factor candidate solar unit $r$.


    \item[$\overline{F}_l$] Capacity of fixed line segment $l$.

    \item[$\overline{F}^{FH}_{l}$] Initial capacity of feeder head transformer.

    \item[$\overline{F}^{FH,upd}_{l}$] Potential upgrade capacity for feeder head

    \item[$\overline{F}^{VR}_{l}$] Capacity of voltage regulator located at line segment $l$.
    
    \item[$\overline{F}^{option}_{l,r}$] Upgrade capacity option $r$ for line segment $l$.

    \item[$\overline{g}^{RTS}_{r}$] Installed capacity of existing rooftop solar unit $r$.

     \item[$\overline{g}^{CS}_{r}$] Capacity of CS project $r$.
   
    \item[$\overline{P}^{in}_h$] Power charging capacity of storage unit $h$.

    \item[$\overline{P}^{out}_h$] Power discharging capacity of storage unit $h$.

    \item[$p^{load}_{n,t}$] Active power load of bus $n$ at time $t$.
    
    \item[$q^{load}_{n,t}$] Reactive power load of bus $n$ at time $t$.


    \item[$R^{option}_{l,r}$] Resistance associated with upgrade option $r$

    \item[$\overline{S}_h$] Duration of storage unit $h$.

    \item[$v^{max}_n$] Maximum voltage of bus $n$.
    
    \item[$v^{min}_n$] Minimum voltage of bus $n$.

    \item[$v^{ref}_n$] Reference feeder voltage.
    

    \item[$X^{option}_{l,r}$] Reactance associated with upgrade option $r$.


    
\end{description}

\subsection*{Decision variables}

\begin{description} [\IEEEsetlabelwidth{500000}\IEEEusemathlabelsep]

    \item[$f^p_{l,t,d}$] Active power flow through line segment $l$.
    
    \item[$f^q_{l,t,d}$] Reactive power flow through line segment $l$.
    
    \item[$\overline{f}^{trf}_{l}$] Capacity of transformer in line segment $l$.
    
    \item[$\overline{f}^{upd}_{l}$] Resulting capacity of candidate line segment $l$.
    
    \item[$g^{CS}_{r,t,d}$] Generation output of candidate solar unit $r$.
    
    \item[$g^{CS,crt}_{r,t,d}$] Curtailed generation output of candidate solar.
    
    \item[$g^{RTS}_{r,t,d}$] Generation output of existing solar unit $r$.
    
    \item[$p^{I,+}_{n,t,d}$] Active power surplus at bus $n$ during time $t$.
    
    \item[$p^{I,-}_{n,t,d}$] Active power deficit at bus $n$ during time $t$.
    
    \item[$p^{in}_{h,t,d}$] Active power charging of storage unit $h$.
    
    \item[$p^{out}_{h,t,d}$] Active power discharging of storage unit $h$.
    
    \item[$p^{subs}_{r,t,d}$] Active power output of substation $r$.
    
    \item[$q^{I,+}_{n,t,d}$] Reactive power surplus at bus $n$.
    
    \item[$q^{I,-}_{n,t,d}$] Reactive power deficit at bus $n$.
    
    \item[$q^{subs}_{r,t,d}$] Reactive power output of substation $r$.
    
    \item[$q^{st,+/-}_{h,t,d}$] Reactive power output of storage unit $h$.
    
    \item[$soc_{h,t}$] State of charge of storage unit $h$ during time $t$.
    
    \item[$soc^{t0}_{h}$] Initial state of charge of storage unit $h$.
    
    \item[$v^{\dagger}_{n,t,d}$] Squared voltage of bus $n$ during time $t$.
    
    \item[$v^{\dagger, mid}_{n,t}$] Auxiliary variable related to the squared voltage of bus $n$ 
    considering an voltage regulator.
    
    \item[$x^{CS,inv}_{r}$] Binary decision variable related to investment in community solar unit $r$. 
    
    \item[$x^{FH,inv}_{l}$] Investment in the feeder head transformer.
    
    \item[$x^{line,inv}_{l,r}$] Investment in option $r$ to upgrade line segment $l$.
    
    \item[$x^{VR,inv}_{l}$] Investment to install a voltage regulator.
    
    
    \item[$x^{st,inv}_{h}$] Investment in storage unit $h$.

\end{description}
\section{Introduction}\label{Introduction}

\IEEEPARstart{I}{n} community solar (CS) programs, a utility or third party - nonprofit or special purpose organization \cite{SunShot} - owns a solar array and sells portions of its power to multiple subscribers, which receive a credit on their electricity bill \cite{KLEIN2021225}. This model has enabled small residential and commercial consumers, including renters, owners of buildings with shaded roofs, or those with limited financial resources, to access clean and renewable energy without the need for personal rooftop panels \cite{NREL_LMI}. 

A particular aspect of these community resources is their deployment at the distribution level. In Europe, for instance, Renewable Energy Communities require members to be in proximity to the renewable projects owned or developed by the community \cite{european_unities_REC}. In the US, an indicator of the local characteristic of CS projects is the capacity of the projects. According to the Community Solar Project Database, 57\% of the projects installed the US have capacities between 0.5-5MW \cite{SS_db}, which indicates that they are often connected to the local distribution grid. Legal and regulatory definitions of CS projects follow similar trends regarding the size of these projects: for example, the state of Maine, in the US, introduced a shared distributed generation procurement process that includes CS assets, developed at the local level, limited to a maximum system size of 5 MW \cite{DL1711}. 

Thus, CS projects must find a balance in terms of capacity, being large enough to create economies of scale and enable efficient asset sharing, yet small enough to be considered local. Unlike behind-the-meter PV, often sized based on consumers' demand and characterized by multiple smaller interconnections at LV (low voltage) level, CS projects require the interconnection of MW scale solar panels at a single point of the distribution grid without a nodal load to offset generation. This unique characteristic poses significant challenges to the distribution grid and has resulted in CS interconnection delays \cite{nyseia_2021} and high interconnection costs \cite{rmi_2018}, which can significantly impact the economic feasibility of these shared solar projects. 

On the other hand, strategically placed distributed PV has the potential to defer the need for distribution grid investments \cite{KEEN2019110902}, presenting a unique opportunity for CS projects. Unlike rooftop solar, CS projects are not tied to a specific feeder location, allowing them to be situated in areas that can either minimize the impact or even benefit the distribution grid. However, in regulatory proceedings, these potential deferral benefits have not been consistently factored into the calculation of CS-related distribution grid costs. 

This paper aims to quantify and discuss distribution grid costs and opportunities associated with CS projects to inform policy and regulatory processes around planning and valuation of these resources. Specifically, we seek to understand: i) the range of costs and associated type of distribution infrastructure investments; ii) and the potential for cost deferral in CS projects, particularly those strategically sited.

\subsection{Literature Review}

In many jurisdictions, the integration of distributed PV systems starts with a hosting capacity analysis, in which utilities quantify the amount of solar that can be interconnected into a feeder without adversely impacting power quality or reliability under existing control and protection systems. Hosting capacity methodologies, such as the one developed in \cite{shayani_hosting}, apply load flow techniques to quantify maximum PV penetration levels within voltage and line capacity limits. The methods can be extended to include operational aspects, such as static VAr compensation and transformers with on-load tap changers (OLTC) control \cite{wang_hosting_capacity}. However, to connect larger capacities of distributed PV projects that surpass the hosting capacity limit, upgrading the existing distribution infrastructure may become necessary. In \cite{Gensollen_2019}, a simulation-based methodology is introduced to assess the cost associated with enhancing grid hosting capacity. Optimization versions of these methods can be also found in stochastic \cite{Santos} and robust \cite{Melgar} forms. However, comprehensive hosting capacity enhancement studies are practically nonexistent, with the exception of the work of \cite{GUPTA2021116010} that provides a nationwide study for hosting capacity enhancement via investments in battery energy storage systems (BESS).

Instead of generically improving hosting capacity, a different body of literature focuses on planning distribution upgrades to integrate concrete PV scenarios. In \cite{osti_1432760}, a set of prototypical feeders are upgraded to integrate distributed PV scenarios, using conventional distribution planning solutions, such as voltage regulators, reconductoring, and transformer upgrades. Non-conventional approaches to address rooftop PV integration, including controls and BESS, are considered in the techno-economic assessment outlined in \cite{Resch}. Both works employ a series of power flows to identify feasible upgrade investments, which are then selected based on least-cost heuristics. Nevertheless, as documented in \cite{HOROWITZ2018420}, these distribution PV integration cost studies do not go beyond a few feeder cases, which may not be comprehensive enough to inform statewide policy and regulatory decisions. A couple of exceptions to this trend include a behind-the-meter PV integration cost study that comprises 75 feeders across three California distribution utilities \cite{dnv_2017} and a recent work that estimates upgrade costs associated with behind-the-meter PV deployment in a distribution network with 170 thousand consumers in Switzerland \cite{GUPTA2021116504}.

Methodologically, an alternative to a sequence of power flows along with heuristics for determining least-cost upgrades is to formalize these problems through distribution grid capacity expansion models, which aim to determine the optimal combination of network upgrades necessary to integrate specific netload profiles (including PV), while ensuring grid constraints. Several distribution capacity expansion models have been proposed, encompassing various investment solutions such as substation/circuit upgrades and distributed generation \cite{Munoz2015}, sizing and placement of BESS \cite{Nick2014}, network reconfiguration \cite{Nick2018}, and OLTC upgrades \cite{IRIA20191147}. These models have also been extended to address objectives beyond cost-minimization, including considerations of reliability \cite{Jooshaki2022} and resilience \cite{Moreira2023}. However, so far, comprehensive studies quantifying distributed PV integration costs have not relied on these capacity expansion models. For example, among the 12 studies detailed in \cite{HOROWITZ2018420}, none incorporated optimization techniques, and only a limited fraction employed load flow calculations.

Finally, empirical data can also be used to determine distribution grid upgrade costs related to distributed PV. In the case of CS, empirical data from the Community Solar Gardens program in Minnesota is one of the few sources of information on these costs \cite{SG_2020, SG_2021}. Examples of empirical analyzes of distributed PV use interconnection and utility investment historical data to estimate distribution system costs \cite{osti_1494635} or deferrals \cite{KEEN2019110902} associated with future solar penetration scenarios. However, these analyzes imply that historical and traditional planning practices will continue in the future, which neglects emerging distribution system planning processes, for example based on non-wire alternatives \cite{rmi_nwa}. 

\subsection{Contributions}
Thus, there is a lack of forward-looking comprehensive studies to assess distribution grid upgrade costs associated with distributed PV and inform the policy and regulatory design of distributed solar programs. Comprehensive studies capable of providing this information are limited to generic hosting capacity enhancements \cite{GUPTA2021116010} or focuses exclusively on behind-the-meter PV \cite{dnv_2017, GUPTA2021116504}. Based on these observation we summarize as follows our contributions in this paper.
\begin{enumerate}
    \item We propose a new methodology to quantify distribution upgrade costs associated with behind-the-meter PV. Different from the simulation approaches used in \cite{dnv_2017, GUPTA2021116504}, this methodology relies on incremental least-cost expansion of the distribution system.
    \item Beyond existing solar integration costs studies, we consider a both traditional and non-traditional investments and deferrals, together with mitigation actions specific to CS projects, including strategic siting and capacity downsizing.
    \item  We apply our methodology to more than 2,500 feeders, considering realistic sizing practices of CS projects in the US, to obtain distribution grid upgrade costs associated with CS integration. We show that our results are realistic by comparing them against 210 real CS integration costs reported in \cite{SG_2020, SG_2021}. We show costs from the consumers' and project developers' perspectives and present policy and regulatory recommendations.
\end{enumerate}

\section{Methodology overview}
\label{sec.MethodologyOverview}

In this section, we provide an overview of our proposed methodology comprising a discussion on (i) the differences between CS and behind-the-meter PV, (ii) the netload scenarios under consideration, (iii) the upgrades and costs, (iv) the methodology workflow, and (v) the granularity of the expansion plan.

\subsection{Community solar vs behind-the-meter PV}

Community solar (CS) can be defined as “{\it any solar project or purchasing program, within a geographic area, in which the benefits of a solar project flow to multiple customers such as individuals, businesses, nonprofits, and other groups}" \cite{CS_program}. In this work, and given CS project sizes in the US (0.5-5MW) \cite{SS_db}, we assume they are connected to the MV portion of the feeder. As CS generation is sold to the grid at a fixed rate and credited to subscribers, developers have an incentive to maximize project capacity. To align with this industry trend, we assume that the CS project's MW capacity matches the minimum daily load (MDL) of the feeder, often used as a cap for distributed PV \cite{NREL_distV}. We also considered that the point of interconnection, depending on the scenario, is either random or strategically located to minimize the system costs.

The proposed cost assessment methodology applies only to this form of distributed PV, (here referred as CS) in which capacity is maximized relative to the feeder and connected at a single (potentially flexible) interconnection point. It is important to stress that these characteristics are different from behind-the-meter (or rooftop) PV, often sized based on load or solar compensation tariffs (e.g. net-metering) and linked to a consumer meter. Thus, behind-the-meter PV is only considered as part of the pre-existing feeder netload scenarios.  

\subsection{Feeder netload scenarios}

The costs to integrate CS into the distribution grid are calculated for different scenarios of feeder loading, reproducing different conditions of the feeder immediately before a CS project is installed. We construct 3 variations of feeder loading and develop the following scenarios: 

\begin{itemize}
    \item Base. This corresponds to the existing conditions loading conditions of the feeder. 
    \item High Rooftop PV. To obtain this scenario, and starting from the ``Base'' netload conditions, we iteratively increase the pre-existing (non-CS) behind-the-meter PV penetration in the feeder. Through a series of power flow simulations, we continue this process until overvoltage or line capacity violations begin to occur. The aim is to simulate conditions close to the feeder hosting capacity, where a substantial amount of PV infrastructure is already in place before the installation of a CS project.
    \item High Load. To achieve this scenario, and starting from the "Base" netload conditions, we incrementally raise the load in the feeder until undervoltage or line capacity violations start to occur. The objective here is to bring the system close to its load capacity limit, thereby simulating scenarios with high loads, such as those resulting from electrification policies. 
\end{itemize}

\subsection{Upgrades and costs}

For each scenario of netload, we calculate the necessary upgrades in the feeder required to maintain the operational limits of voltage and power in 2 situations: without and with CS. The type of upgrades include reconductoring of overhead and underground lines, replacement of transformers, installation of voltage regulators, and utility-owned BESS. For the case of "with CS", two additional mitigation strategies related to the CS project are considered: a) downsizing the CS capacity; b) strategically siting the project.

The combination of upgrade and mitigation strategies for each scenario are determined based on a least-cost optimization model, presented in section \ref{sec.MathematicalFormulation}. The overall integration costs of community solar projects ($C^{itgr}_{CS}$) correspond to the difference between the optimal costs with ($C_{CS}$) and without ($C_{\overline{CS}}$) the project in the system, resulting from the objective function described in section \ref{sec.MathematicalFormulation}.

\begin{align}
    C^{itgr}_{CS} = C_{CS} - C_{\overline{CS}} \label{eqn:incremental_cost}
\end{align}

Thus, $C^{itgr}_{CS}$ captures the incremental costs of distribution grid upgrades that are due to the presence of CS in the system. A positive incremental cost means that CS will require overall additional infrastructural upgrades. Conversely, a negative incremental cost signifies that the CS project is anticipated to result in distribution grid investment deferrals.

\subsection{Methodology workflow}

The first step of the cost assessment is to determine the 3 netload scenarios presented above. This is accomplished by executing homothetic variations of netload, coupled with three-phase power flow evaluations for each critical hours. During peak load hours, the substation transformer tap is raised, and capacity and undervoltage limits are evaluated. Conversely, during maximum PV production hours, the substation transformer tap is reduced, and overvoltage limits are assessed.

The three-phase power flow results and the violation risk assessment, including lines and nodes near technical limits, help inform the selection of potential candidate locations for upgrades: line and transformers near their capacity are considering for reconductoring and capacity upgrades, respectively; nodes at risk of voltage violations will have their neighboring nodes, both upward and downward, identified as candidates for the placement of voltage regulators. Besides this set of candidates, 3 locations - beginning (close to the substation), middle and end of the feeder - are selected for potential placement of BESS. In the "with CS" cases, the same 3 locations are considered for siting CS projects.

Considering the netload scenario and the resulting set of candidates, we run a least-cost optimal distribution grid expansion model, whose formulation is presented in section \ref{sec.MathematicalFormulation}. At the end of the run, the slack variables of the model are checked and if violations persist the list of candidates is updated using the process described above. This iterative cycle continues until all slack variables reach zero. The process ends with a verification of the system security using a three-phase power flow analysis. 

\begin{figure}[!ht]
       \centering
      \includegraphics[width=0.45 \textwidth]{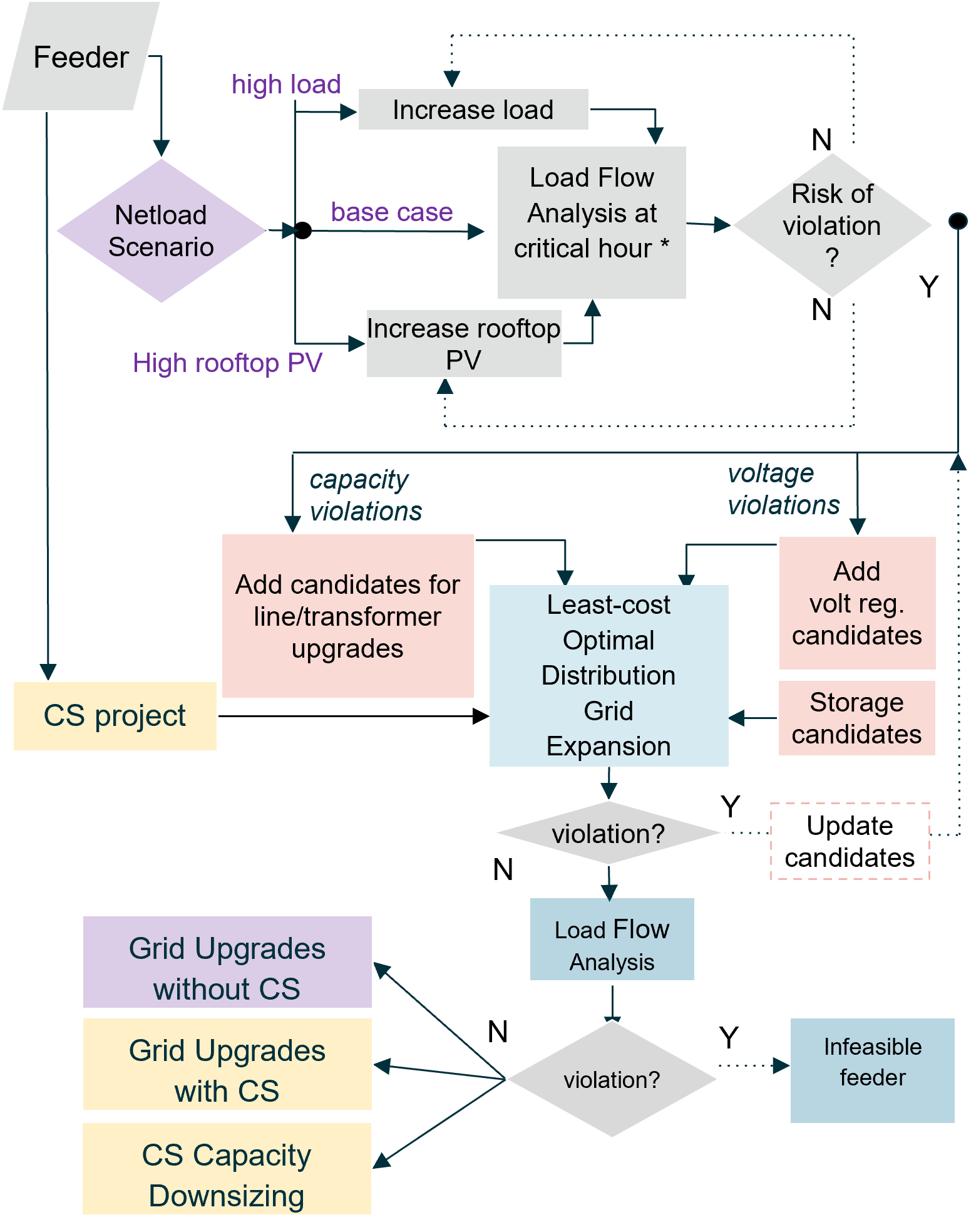}
      \caption{Flowchart of the methodology.}
      \label{Fig: methodology_overview}
\end{figure}

Figure \ref{Fig: methodology_overview} presents an overview of workflow used in the CS distribution grid cost assessment.

\subsection{Granularity of the expansion plan}

The cost assessment is performed at a feeder level. Each feeder expansion plan is run with an hourly resolution and considering three typical days in a year: the day of the peak, the day with the maximum solar production, and an average day. The first two days, each assigned a weight of 1, ensure the planning solution's feasibility under extreme netload conditions. The third day, with a weight of 363, captures the average daily economic of operations of the system, including BESS set-points and CS generation curtailment, on the planning decisions.

Given that CS projects and BESS are linked to the MV three-phase section of the feeder, their presence does not impact investment costs in downward single-phase laterals and the low voltage portion of the circuit. Consequently, according to equation (\ref{eqn:incremental_cost}), investments in these sections of the circuit do not affect CS integration costs. Therefore, we aggregate these circuits and their netloads at the MV three-phase nodes. This not only reduces the number of nodes in the expansion plan but also allows the use of balanced approximations of the feeder power flow constraints in the expansion plan. It is important to note that, despite these approximations, the security of the planning solutions is validated against a full three-phase unbalance power flow.

\section{Least-cost optimal distribution grid expansion formulation}\label{sec.MathematicalFormulation}

The mathematical model formulated in this section, aims at obtaining the most cost-effective portfolio of investments in distribution grids to integrate CS. The upgrades include reinforcement of line segments (reconductoring), voltage regulators, BESS and transformer upgrades. Next, we provide a brief setup describing the relationships between the different sets of line segments and, then, describe the formulation in detail.

\textbf{Setup.} Consider ${\cal L} = {\cal L}^{C} \cup {\cal L}^{F} \cup {\cal L}^{FH} \cup {\cal L}^{VR}$, ${\cal L}^{C} \cap {\cal L}^{F} \cap {\cal L}^{FH} \cap {\cal L}^{VR} = \emptyset$, ${\cal L}^{C} \cap {\cal L}^{F}= \emptyset$, ${\cal L}^{C} \cap {\cal L}^{FH}= \emptyset$, ${\cal L}^{C} \cap {\cal L}^{VR}= \emptyset$, ${\cal L}^{F} \cap {\cal L}^{FH}= \emptyset$, ${\cal L}^{F} \cap {\cal L}^{VR}= \emptyset$, ${\cal L}^{FH} \cap {\cal L}^{VR}= \emptyset$. Also, ${\cal L}^{VR} = {\cal L}^{VR,E} \cup {\cal L}^{VR,C}$, $v^{\dagger}_{nt} = v^{2}_{nt}$. For lines with segments with tap changers, we have $v^{\dagger,mid}_{to(l),t} = v^{\dagger}_{to(l),t,d} \varphi^2_l$, where $\varphi^2_l$ and $to(l)$ are the turns ratio and $to(l)$ is the receiving bus of the transformer in $l$, respectively.


\subsection{Objective} \label{ObjFunc_Model_InvGenBat}

\begin{align}
    &\underset{{\substack{ f^p_{l,t,d}, f^q_{l,t,d}, \overline{f}^{trf}_{l}, \overline{f}^{upd}_{l}, g^{CS}_{r,t,d}, g^{CS,crt}_{r,t,d},\\ g^{RTS}_{r,t,d}, p^{I,+}_{n,t,d}, p^{I,-}_{n,t,d}, p^{in}_{h,t,d}, p^{out}_{h,t,d}, p^{subs}_{r,t,d},\\ q^{I,+}_{n,t,d}, q^{I,-}_{n,t,d}, q^{subs}_{r,t,d}, q^{st,+/-}_{h,t,d}, soc_{h,t,d}, \\ soc^{t0}_{h}, v^{\dagger}_{n,t,d}, v^{\dagger, mid}_{n,t,d}, x^{CS,inv}_{r}, x^{FH,inv}_{l,r},\\ x^{line,inv}_{l,r}, x^{VR,inv}_{l}, x^{st,bin,inv}_{h}, x^{st,inv}_{h} }}} {\text{Minimize}}  \hspace{-0.2cm} \sum_{h \in H^{C}} C^{st,inv}_{h} x^{st,inv}_{h} \notag \\
    &\hspace{20pt} + \sum_{l \in {\cal L}^{VR,C}} C^{VR,inv}_l x^{VR,inv}_l \notag \\
    &\hspace{20pt} + \sum_{l \in {\cal L}^C} \sum_{r \in R^{line,options}_l} C^{line, inv}_{l,r} x^{line,inv}_{l,r}  \overline{F}^{option}_{l,r} \notag\\
    &\hspace{20pt} + 
    {\color{black} \sum_{l \in {\cal L}^{FH}} C^{FH,inv}_{l} } x^{FH,inv}_{l} 
    \notag\\
    &\hspace{20pt} + \sum_{d \in {\cal D}} \omega_d \sum_{t \in T}\sum_{r \in R^{Solar,C}} C^{CS,crt}_t   g^{CS,crt}_{r,t,d} \label{Model_v1_1} 
    \notag\\ 
    &\hspace{20pt} + \sum_{d \in {\cal D}} \omega_d \sum_{t \in T}\sum_{n \in N} C^{I}  (p^{I,+}_{n,t,d} + p^{I,-}_{n,t,d} + q^{I,+}_{n,t,d} + q^{I,-}_{n,t,d}) 
\end{align}

The first three terms of the objective function in \eqref{Model_v1_1} comprise costs related to storage placement, voltage regulator installation and reconductoring. Since the analysis is performed at the feeder level, we apportion transformer upgrades relative to the feeder. The fourth term of the objective function models the cost of these "feeder head transformer reinforcement". The  last two terms model cost of CS generation curtailment and the cost of load shedding slack variables.

\subsection{Nodal balance and reference voltage}\label{Balance_Model_InvGenBat}
\begin{align}
    & \sum_{r \in R^{Subs}_{n}} p^{subs}_{r,t,d} + \sum_{l \in {\cal L}^{to}_{n}} f^p_{l,t,d} - \sum_{l \in {\cal L}^{fr}_{n}} f^p_{l,t,d} + \sum_{r \in R^{Solar,C}_{n}} g^{CS}_{r,t,d} \notag \\ 
    & \hspace{5pt} + \sum_{r \in R^{Solar,E}_{n}} g^{RTS}_{r,t,d} + \sum_{h \in H_{n}} p^{out}_{h,t,d} - \sum_{h \in H_{n}} p^{in}_{h,t,d} - p^{load}_{n, t, d} \notag \\ 
    & \hspace{5pt} - \sum_{l \in {\cal L}^{to}_{n}} \frac{\beta^p_{l}}{2} \Biggl [ \sum_{n \in N} p^{load}_{n, t, d}- \sum_{r \in R^{Solar,E}} g^{RTS}_{r,t,d} - \sum_{r \in R^{Solar,C}_{n}} g^{CS}_{r,t,d} \Biggr ] \notag \\ 
    & \hspace{5pt}- \sum_{l \in {\cal L}^{fr}_{n}} \frac{\beta^p_{l}}{2} \Biggl [ \sum_{n \in N} p^{load}_{n, t, d}- \sum_{r \in R^{Solar,E}} g^{RTS}_{r,t,d} - \sum_{r \in R^{Solar,C}_{n}} g^{CS}_{r,t,d} \Biggr ]  \notag \\ 
    & \hspace{48pt} + p^{I,-}_{n,t,d} - p^{I,+}_{n,t,d} = 0;  \forall n \in N, t \in T, d \in {\cal D} \label{Model_v1_2} \\
    & \sum_{r \in R^{Subs}_{n}} q^{subs}_{r,t,d} + \sum_{l \in {\cal L}^{to}_{n}} f^q_{l,t,d} - \sum_{l \in {\cal L}^{fr}_{n}} f^q_{l,t,d} + \sum_{h \in H_{n}} q^{st,+/-}_{h,t,d} \notag \\ 
    & \hspace{5pt} - q^{load}_{n, t,d} - \sum_{l \in {\cal L}^{to}_{n}} \frac{\beta^q_{l}}{2} \sum_{n \in N} q^{load}_{n, t,d} - \sum_{l \in {\cal L}^{fr}_{n}} \frac{\beta^q_{l}}{2}  \sum_{n \in N} q^{load}_{n, t,d} \notag \\ 
    & \hspace{48pt} + q^{I,-}_{n,t,d} - q^{I,+}_{n,t,d} = 0;   \forall n \in N, t \in T, d \in {\cal D} \label{Model_v1_3} \\
    & v^{\dagger}_{fr(l),t,d} = (v^{ref})^2; \forall l \in {\cal L}^{FH}, t \in T, d \in {\cal D} \label{Model_v1_4} \\
    & v^{{min}^2}_n \leq  v^{\dagger}_{n,t,d} \leq v^{{max}^2}_n; \forall n \in N, t \in T, d \in {\cal D} \label{Model_v1_5} 
\end{align}

Constraints \eqref{Model_v1_2} enforce active nodal balance for buses in the system considering potential local injection, storage devices, power flows in and out, and the load of the bus, whereas constraints \eqref{Model_v1_3} play an analogous role for reactive balance. Voltage reference for feeder heads and voltage limits for all buses are imposed in \eqref{Model_v1_4} and \eqref{Model_v1_5}, respectively.

\subsection{Fixed lines segments}


\begin{align}
    &  - \overline{F}_l \leq  f^p_{l,t,d} \leq \overline{F}_l; \forall l \in {\cal L}^{F}, t \in T, d \in {\cal D} \label{Model_v1_6}  \\
    &  - \overline{F}_l \leq  f^q_{l,t,d} \leq \overline{F}_l; \forall l \in {\cal L}^{F}, t \in T, d \in {\cal D} \label{Model_v1_7} \\
    &  f^q_{l,t,d} \leq \overline{f}^q(f^p_{l,t,d},\overline{F}_l,e); \forall l \in {\cal L}^{F}, t \in T, d \in {\cal D}, \notag \\
    & \hspace{165pt} e \in \{1,\dots,4\} \label{Model_v1_8}\\
    &  -f^q_{l,t,d} \leq \overline{f}^q(f^p_{l,t,d},\overline{F}_l,e); \forall l \in {\cal L}^{F}, t \in T, d \in {\cal D}, \notag \\
    & \hspace{165pt} e \in \{1,\dots,4\} \label{Model_v1_9} \\
    & v^{\dagger}_{to(l),t,d} - \biggl [ v^{\dagger}_{fr(l),t,d} - {2} (R_{l} f^p_{l,t,d} + X_{l} f^q_{l,t,d}) \biggr ]  = 0; \notag \\
    & \hspace{130pt}\forall l \in {\cal L}^{F}, t \in T, d \in {\cal D}  \label{Model_v1_10}
\end{align}

\noindent where $\overline{f}^q(f^p_{l,t,d},\overline{F}_l,e) = cotan \biggl (  \biggl (  \frac{1}{2} - e  \biggr ) \frac{\pi}{4} \biggr ) \biggl ( f^p_{l,t,d} - cos \biggl (e \frac{\pi}{4} \biggr )  \overline{F}_l  \biggr ) + sin \biggl( e \frac{\pi}{4} \biggr) \overline{F}_l$.

Constraints \eqref{Model_v1_6} and \eqref{Model_v1_7} enforce bounds for active and reactive power flows, respectively. Following \cite{Mashayekh2018}, we also impose limits on active and reactive power flows together in \eqref{Model_v1_8} and \eqref{Model_v1_9} based on a linear approximation of the apparent power flow and describe voltage drop through expression \eqref{Model_v1_10}.   


\subsection{Candidates line segments for reconductoring}

\begin{align}
  & \overline{f}^{upd}_l = \sum_{r \in R^{line,options}_l} x^{line,inv}_{l,r} \overline{F}^{option}_{lr}; \forall l \in {\cal L}^C  \label{Model_v1_11}\\
  &\sum_{r \in R^{line,options}_l} x^{line,inv}_{l,r} = 1; \forall l \in {\cal L}^C  \label{Model_v1_12}\\
  & x^{line,inv}_{l,r} \in \{0,1\}; \forall l \in {\cal L}^C, r \in R^{line,options}_l  \label{Model_v1_13}\\
  & - \overline{f}^{upd}_l \leq  f^p_{l,t,d} \leq \overline{f}^{upd}_l; \forall l \in {\cal L}^C, t \in T, d \in {\cal D}  \label{Model_v1_14}\\
  & - \overline{f}^{upd}_l \leq  f^q_{l,t,d} \leq \overline{f}^{upd}_l; \forall l \in {\cal L}^C, t \in T, d \in {\cal D} \label{Model_v1_15}\\
  & f^q_{l,t,d} \leq \overline{f}^q(f^p_{l,t,d},\overline{f}^{upd}_l,e);\forall l \in {\cal L}^C, t \in T, d \in {\cal D}, \notag\\
  & \hspace{165pt} e \in \{1,\dots,4\} \label{Model_v1_16}\\
  & -f^q_{l,t,d} \leq \overline{f}^q(f^p_{l,t,d},\overline{f}^{upd}_l,e); \forall l \in {\cal L}^C, t \in T, d \in {\cal D}, \notag\\
  & \hspace{165pt} e \in \{1,\dots,4\}   \label{Model_v1_17}\\
  & - (1-x^{line,inv}_{lr}) M \leq v^{\dagger}_{to(l),t,d} - \biggl [ v^{\dagger}_{fr(l),t,d}  \notag\\
  & \hspace{2pt} - {2} (R^{option}_{lr} f^p_{l,t,d} + X^{option}_{lr} f^q_{l,t,d}) \biggr ]  \leq(1 -x^{line,inv}_{lr}) M; \notag\\
  & \hspace{60pt}  \forall l \in {\cal L}^C, t \in T, d \in {\cal D}, r \in R^{line,options}_l \label{Model_v1_18}
\end{align}

Expressions \eqref{Model_v1_11}--\eqref{Model_v1_13} model the potential selection of an upgrade option for each line segment $l \in {\cal L}^C$. It is worth mentioning that one of the options is to keep line segment $l$ unaltered, whose respective cost $C^{line,inv}_{l,r}$ is null. Depending on the values assumed by decision variables $x^{line,inv}_{l,r}$, the corresponding power flow limits will be considered in \eqref{Model_v1_14}--\eqref{Model_v1_17} and the resistances and reactances will be taken into account for voltage drop in \eqref{Model_v1_18}.   

\subsection{Tap changing}
 
\begin{align}
    &\frac{v^{\dagger, mid}_{to(l), t, d}}{(\varphi^{max}_l)^2} \leq v^{\dagger}_{to(l), t, d} \leq \frac{v^{\dagger, mid}_{to(l), t, d}}{(\varphi^{min}_l)^2}; \forall l \in {\cal L}^{VR} \cup {\cal L}^{FH},\notag\\
    &\hspace{170pt} t \in T, d \in {\cal D} \label{Model_v1_19} \\
    & - x^{VR,inv}_l M  \leq v^{\dagger, mid}_{to(l),t,d} - v^{\dagger}_{to(l), t, d} \leq x^{VR,inv}_l M; \notag\\
    &\hspace{115pt} \forall l \in {\cal L}^{VR,C}, t \in T, d \in {\cal D}\label{Model_v1_20}\\
    & x^{VR,inv}_l \in \{0,1\}; \forall l \in {\cal L}^{VR,C} \label{Model_v1_21}
\end{align}

We consider the possibility of tap changing in voltage regulator and OLTC transformers located at the feeder head. In this context, expression \eqref{Model_v1_19} imposes limits of tap changing. In addition expressions \eqref{Model_v1_20} and \eqref{Model_v1_21} model the logic behind enabling tap changing for the candidate voltage regulators that are selected for installation. 

\subsection{Feeder head and voltage regulators}

\begin{align}
    & \overline{f}^{trf}_l = \overline{F}^{FH}_l + x^{FH,inv}_l \overline{F}^{FH,upd}_l; \forall l \in {\cal L}^{FH} \label{Model_v1_22}\\
    & \overline{f}^{trf}_l = \overline{F}^{VR}_l; \forall l \in {\cal L}^{VR} \label{Model_v1_23}\\
    & x^{FH,inv}_l \in \{0,1\}; \forall l \in {\cal L}^{FH} \label{Model_v1_24} \\
    &  - \overline{f}^{trf}_l \leq  f^p_{l,t,d} \leq \overline{f}^{trf}_l; \forall l \in {\cal L}^{VR} \cup {\cal L}^{FH}, t \in T, \notag\\
    & \hspace{190pt} d \in {\cal D} \label{Model_v1_25}\\
    &  - \overline{f}^{trf}_l \leq  f^q_{l,t,d} \leq \overline{f}^{trf}_l; \forall l \in {\cal L}^{VR} \cup {\cal L}^{FH}, t \in T, \notag\\
    & \hspace{190pt} d \in {\cal D} \label{Model_v1_26}\\
    &  f^q_{l,t,d} \leq \overline{f}^q(f^p_{l,t,d},\overline{f}^{trf}_l,e); \forall l \in {\cal L}^{VR} \cup {\cal L}^{FH}, t \in T, \notag\\
    & \hspace{135pt} d \in {\cal D}, e \in \{1,\dots,4\} \label{Model_v1_27}\\
    &  -f^q_{l,t,d} \leq \overline{f}^q(f^p_{l,t,d},\overline{f}^{trf}_l,e); \forall l \in {\cal L}^{VR} \cup {\cal L}^{FH}, t \in T, \notag\\
    & \hspace{135pt} d \in {\cal D}, e \in \{1,\dots,4\}  \label{Model_v1_28}\\
    & v^{\dagger, mid}_{to(l),t,d} - \biggl [ v^{\dagger}_{fr(l),t,d} - {2} (R_{l} f^p_{l,t,d} + X_{l} f^q_{l,t,d}) \biggr ]   = 0;\notag\\
    & \hspace{90pt} \forall l \in {\cal L}^{VR} \cup {\cal L}^{FH}, t \in T, d \in {\cal D} \label{Model_v1_29} 
\end{align}

Expressions \eqref{Model_v1_22} and \eqref{Model_v1_23} determine the capacities of feeder head transformers, which can be upgraded according to \eqref{Model_v1_24}, and voltage regulators, respectively. Based on these capacities, power flow limits are enforced in \eqref{Model_v1_25}--\eqref{Model_v1_28}. Moreover, constraints \eqref{Model_v1_29} represent voltage drop. 

\subsection{Solar generation} \label{Solar_Model_InvGenBat}

\begin{align}
    & g^{RTS}_{r,t,d} =  \overline{g}^{RTS}_{r}  f^{RTS}_{r,t,d}; \forall r \in R^{Solar,E}, t \in T, d \in {\cal D}  \label{Model_v1_30}  \\
    & g^{CS}_{r,t,d} + g^{CS,crt}_{r,t,d} =  x^{CS,inv}_{r} \overline{g}^{CS}_{r} f^{CS}_{r,t,d}; \forall r \in R^{Solar,C}, \notag\\
    & \hspace{170pt} t \in T, d \in {\cal D}  \label{Model_v1_31}\\
    & \sum_{r\in R^{Solar,C}} x^{CS,inv}_r  =  1 \label{Model_v1_32} \\
    & x^{CS,inv}_r  \in \{0,1\}; \forall r \in R^{Solar,C} \label{Model_v1_32_a} 
\end{align}

For rooftop solar, expression \eqref{Model_v1_30} determines the available existing generation to be absorbed by the system during each time $t$. Analogously, constraints \eqref{Model_v1_31} inform the available output for CS project. The location of CS is modeled by the binary variable $x^{CS,inv}_r$ in \eqref{Model_v1_32} and \eqref{Model_v1_32_a}, when CS is allowed to be strategically sited. Otherwise $x^{CS,inv}_r$ becomes a parameter.

\subsection{Operation of storage devices} \label{Storage_Model_InvGenBat}

\begin{align}
    & soc^{t0}_{h,d} =  soc_{h, t^{last}_d,d};\forall h \in H, d \in {\cal D}  \label{Model_v1_33}  \\
    & soc_{h, t,d} =  soc^{t0}_{h,d} + \eta  p^{in}_{h,t,d} - p^{out}_{h,t,d}; \forall h \in H, d \in {\cal D}, \notag\\
    & \hspace{190pt} t = t^{ini}_d  \label{Model_v1_34}  \\
    & soc_{h, t,d} =  soc_{h, t-1,d} + \eta p^{in}_{h,t,d} - p^{out}_{h,t,d}; \forall h \in H, \notag\\
    & \hspace{135pt} t \in T, d \in {\cal D} | t \neq t^{ini}_d   \label{Model_v1_35}  \\
    & 0 \leq soc_{h, t, d} \leq  \overline{S}_h \overline{P}^{in}_h; \forall h \in H^{E}, t \in T, d \in {\cal D}  \label{Model_v1_36}  \\
    & 0 \leq soc_{h, t, d} \leq  \overline{S}_h x^{st,inv}_h \overline{P}^{in}_h; \forall h \in H^{C}, t \in T, d \in {\cal D}  \label{Model_v1_37}  \\
    & 0 \leq p^{in}_{h, t, d} \leq \overline{P}^{in}_h; \forall h \in H^{E}, t \in T, d \in {\cal D}  \label{Model_v1_38}  \\
    & 0 \leq p^{in}_{h, t, d} \leq  x^{st,inv}_h \overline{P}^{in}_h; \forall h \in H^{C}, t \in T, d \in {\cal D}  \label{Model_v1_39} \\
    & 0 \leq p^{out}_{h, t, d} \leq \overline{P}^{out}_h; \forall h \in H^{E}, t \in T, d \in {\cal D}  \label{Model_v1_40}  \\
    & 0 \leq p^{out}_{h, t, d} \leq  x^{st,inv}_h \overline{P}^{out}_h; \forall h \in H^{C}, t \in T, d \in {\cal D}  \label{Model_v1_41} \\
    & - \psi \overline{P}^{in}_{h} \leq q^{st,+/-}_{h,t,d}  \leq \psi \overline{P}^{in}_{h}; \forall h \in H^{E}, t \in T, d \in {\cal D} \label{Model_v1_42}\\
    & - \psi x^{st,inv}_h \overline{P}^{in}_{h}  \leq q^{st,+/-}_{h,t,d} \leq \psi x^{st,inv}_h \overline{P}^{in}_{h}; \forall h \in H^{C}, \notag\\
    & \hspace{170pt}t \in T, d \in {\cal D} \label{Model_v1_43} \\
    & 0 \leq x^{st,inv}_{h} \leq x^{st,bin,inv}_{h} \overline{x}^{st,inv}_{h}; \forall h \in H^{C} \label{Model_v1_44}\\
    & \sum_{h \in H^C} x^{st,bin,inv}_h = 1 \label{Model_v1_45}\\
    & x^{st,bin,inv}_h \in \{0, 1\}; \forall h \in H^C \label{Model_v1_46}
\end{align}

Constraints \eqref{Model_v1_33}--\eqref{Model_v1_35} model state of charge updates throughout time periods. Constraints \eqref{Model_v1_36} and \eqref{Model_v1_37} impose bounds on state of charge for existing and candidate storage devices, respectively. In addition, active power charging and discharging limits are enforced by constraints \eqref{Model_v1_38}--\eqref{Model_v1_41} also for existing and candidate storage devices, whereas \eqref{Model_v1_42} and \eqref{Model_v1_43} play a similar role regarding reactive power. Finally, expressions \eqref{Model_v1_44}--\eqref{Model_v1_46} model selection and sizing of candidate storage devices.

\section{~Results}


For this study we used the NREL's SMART-DS dataset \cite{smart_ds}, which provides high-quality standardized distribution network models that are geolocated and include a full representation of electrical parameters from the low-voltage customer connections, loads and PV productivity profiles through the MV primary and up to sub-transmission, including substation transformers. The SMART-DS covers 2711 feeders representing rural and urban feeders from entire large geographic regions, including Austin, Texas; Greensboro, North Carolina and the extended San Francisco Bay Area, California. 
We were able to complete the analysis on 95-98\% (depending on the scenario) of those feeders. A small fraction was not included in the analysis, since either an initial power flow did not converge, or no solution existed for the least cost optimal distribution grid expansion model.

Unit costs for distribution system upgrades were taken from the NREL's Distribution Grid Integration Unit Cost Database \cite{cost_db}, which contain a set of real distribution component cost that were gathered exactly with the intention of estimating distribution grid integration costs. We considered BESS costs of \$634/kW for a 2-hour battery with 15 year lifetime and the set of transformer listed in Table~\ref{tab:transformer_costs}. Regarding voltage regulators and reconductoring, given the significantly larger cost reported in California (CA) in comparison with other States,  we considered specific costs for this state. Voltage regulator costs were considered to be \$221.7k in CA and \$38.5k in other states. The reconductoring costs, with values varying with the conductor type and application, are presented in Table~\ref{tab:reconductoring_costs}. For each conductor, the impedance and ampacity parameters were taken from \cite{Kersting_book}. Finally, the costs of CS curtailment were considering equal to the average hour wholesale energy prices for each region, taken from the Cambium dataset \cite{osti_1915250}.

\begin{table}[htbp]
 \footnotesize
  \centering
  {\color{black}
  \caption{Transformer parameters}
    \begin{tabular}{c|c|c|c}
    \hline
    \multicolumn{1}{r}{} & Costs (k\$) & Capacity (MVA) & Lifetime (yr) \\
    \hline
Transformer A &	250 &	11 &	30 \\
Transformer B &	600 &	15 &	30 \\
Transformer C &	947	&   20 &	30 \\
Transformer D &	1400 &	35 &	30 \\
Transformer E &	1900 &	40 &	30 \\
Transformer F &	5000 &	100 &  30 \\
Transformer G &	11000 &	500	& 30 \\
    \hline
    \end{tabular}%
  \label{tab:transformer_costs}%
  }
\end{table}%

\begin{table}[htbp]
 \footnotesize
  \centering
  {\color{black}
  \caption{Overhead (OH) and Underground (UG) reconductoring costs per conductor type (k\$/mile) for rural (R) and urban (U) lines}
    \begin{tabular}{c|c|c|c|c|c|c}
    \hline
    \multicolumn{3}{c}{} California & \multicolumn{3}{c}{}  Non-California \\
    \hline
    conductor &	R-OH  & U-OH &	U-UG &	R-OH & U-OH & U-UG \\ 
    \hline
     ACSR \#4 &	892  & 1,510 &	1,167 &	644 & 1,088 & 174 \\ 
    ACSR \#2 &	1,133 &	1,918 &	1,482 &	818 & 1,381 & 221 \\
    ACSR 1/0 &	1,704 &	2,884 &	2,229 &	1,230 & 2,077 &	333 \\
    ACSR 3/0 &	2,597 &	4,394 &	3,396 &	1,875 &	3,165 & 507 \\
    ACSR 4/0 &	3,248 &	5,497 &	4,247 & 2,345 & 3,959 &	634 \\
    ACSR 336.4 & 5,033 & 8,517 &	6,581 & 3,633 &	6,135 &	983 \\
    ACSR 477 &	6,978 & 11,809 &9,125& 5,037  &	8,506 &	1,363 \\ 
    \hline
    \end{tabular}%
  \label{tab:reconductoring_costs}%
  }
\end{table}%

\begin{figure}[!ht]
    \centering
    \subfigure[Cost/deferral distribution]{
        \centering
        \includegraphics[width=0.25 \textwidth]{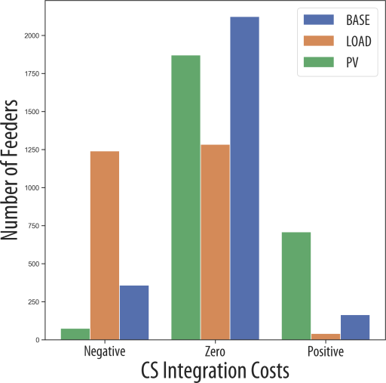}
        \label{fig:cost_distribution}
    }
        \subfigure[Cost comparison]{
        \centering
        \includegraphics[width=0.20 \textwidth]{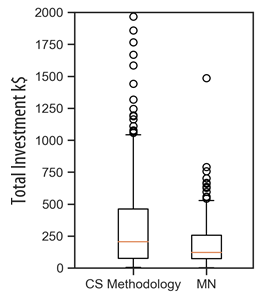}
        \label{Fig: costs_comparison}
    }
      \caption{(a) distribution of CS integration costs obtained across feeders for different netload scenarios, Base, High Rooftop PV, High Load. (b) comparison of positive cost distribution with real CS integration costs reported in MN \cite{SG_2020,SG_2021}.}
\end{figure}

\subsection{Distribution upgrade costs and deferrals}

We categorize distribution upgrade costs associated with CS interconnection into three groups, depending on the value of $C^{itgr}_{CS}$ (see equation \ref{eqn:incremental_cost}): a) \textbf{negative} costs occur when a feeder with an operational CS project defers or avoids grid upgrades that would otherwise be necessary; b) \textbf{zero} costs apply when feeders can accommodate a CS project without requiring additional upgrades compared to the same feeder without CS: c) \textbf{positive} costs arise when CS projects need grid upgrades with a higher cost than the upgrades required (if any) on the same feeder without the CS project. As seen in Figure \ref{fig:cost_distribution}, under "Base" scenario, most feeders do not demand extra distribution grid infrastructure investments for CS projects. When feeders become highly loaded, approximately 50\% experience negative distribution grid infrastructure costs (investment deferral). For feeders with a high penetration of rooftop PV, CS projects necessitate new grid investments in about 30\% of cases. 

Regarding actual values when costs are positive, Figure \ref{Fig: costs_comparison} shows that the average distribution feeder upgrade cost is \$188k in our methodology (labeled as "CS Methodology"). However, in a small number of extreme cases, upgrade costs may reach \$2M (e.g., when a substation transformer upgrade is required). As shown in the same figure \ref{Fig: costs_comparison}, this cost distribution closely aligns with the real upgrade costs of 210 CS projects reported under the Solar Gardens program in Minnesota \cite{SG_2020, SG_2021}. This similarity in cost ranges indicates that our methodology realistically captures the key drivers of CS integration costs.

\begin{figure}[!t]
    \centering
    \subfigure[Costs (\$/kW of CS)]{
        \centering
        \includegraphics[width=0.22 \textwidth]{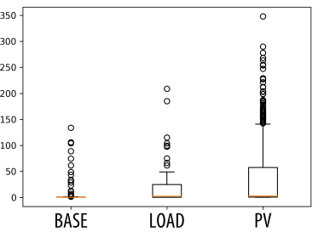}
        \label{fig:cost_cs}
    }
        \subfigure[Deferrals (\$/kW of CS)]{
        \centering
        \includegraphics[width=0.22 \textwidth]{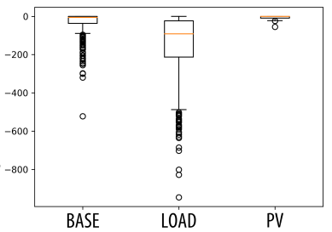}
        \label{Fig: deferrals_cs}
    }
    \\
    \subfigure[Cost (~\cent/kWh of load)]{
        \centering
        \includegraphics[width=0.22 \textwidth]{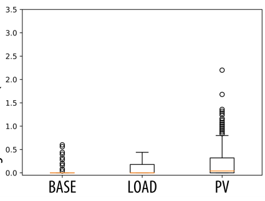}
        \label{fig:cost_cons}
    }
        \subfigure[Deferrals (~\cent/kWh of load)]{
        \centering
        \includegraphics[width=0.22 \textwidth]{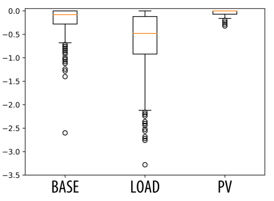}
        \label{Fig: deferrals_cons}
    }
    
      \caption{Distribution of costs and deferrals relative to the size of CS projects (a,b) and relative to the energy demand (c,d).}
\end{figure}

\begin{figure}[!t]
    \centering
        \subfigure[Incremental investments by type]{
        \centering
        \includegraphics[width=0.22 \textwidth]{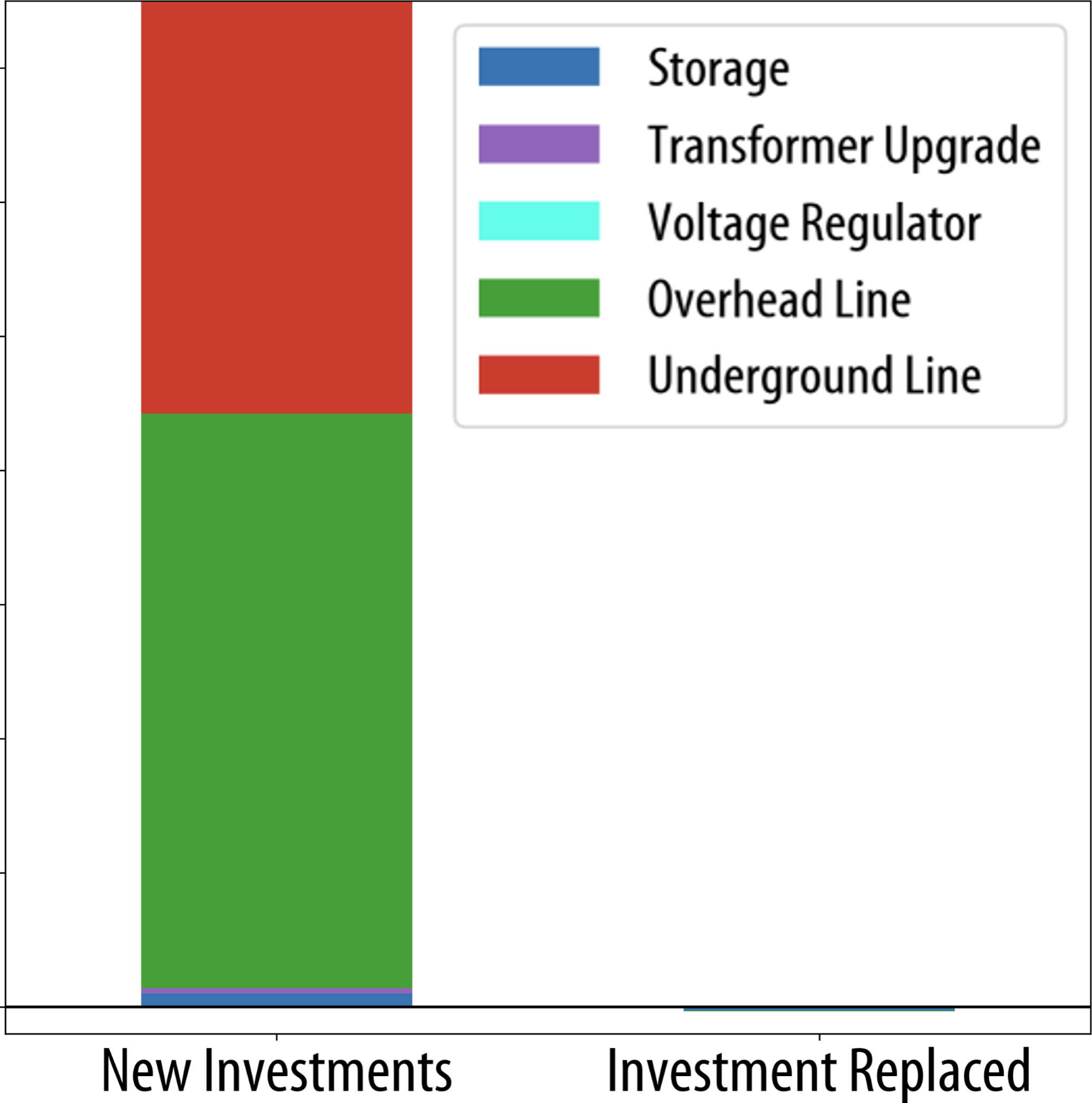}
        \label{fig:incremental_investments}
    }
        \subfigure[Investment type per feeder]{
        \centering
        \includegraphics[width=0.22 \textwidth]{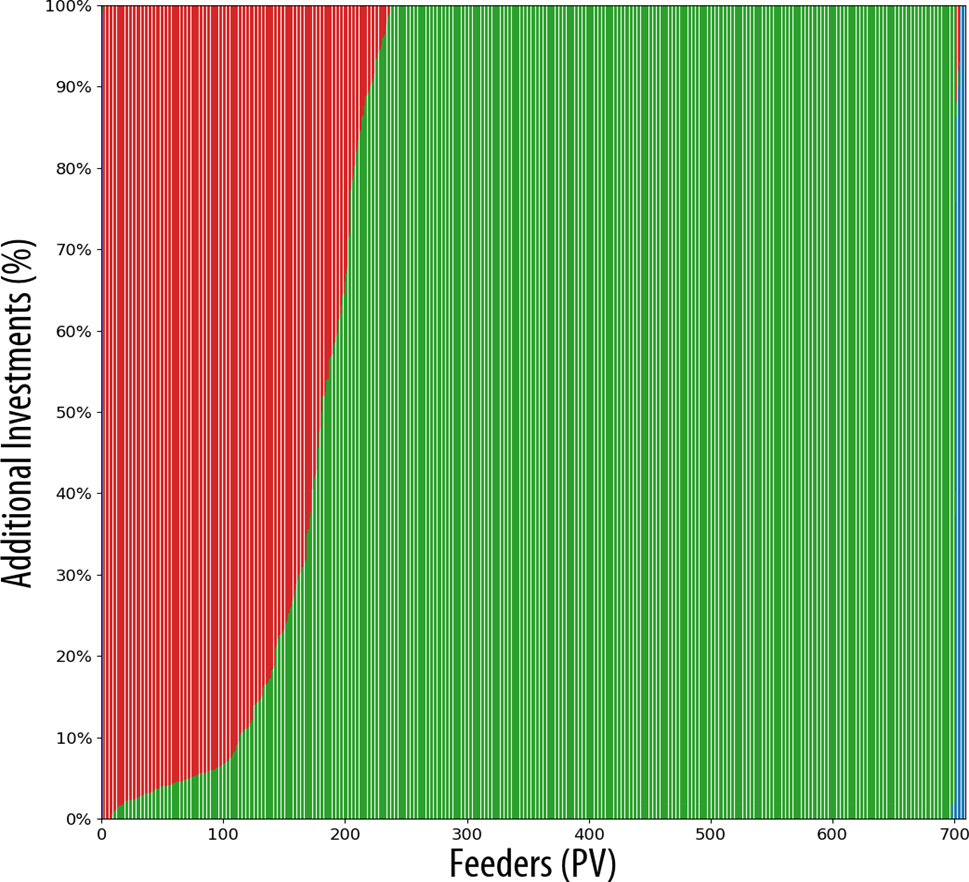}
        \label{Fig:investments_feeder}
    }
    \\
    \subfigure[Incremental deferrals by type]{
        \centering
        \includegraphics[width=0.22 \textwidth]{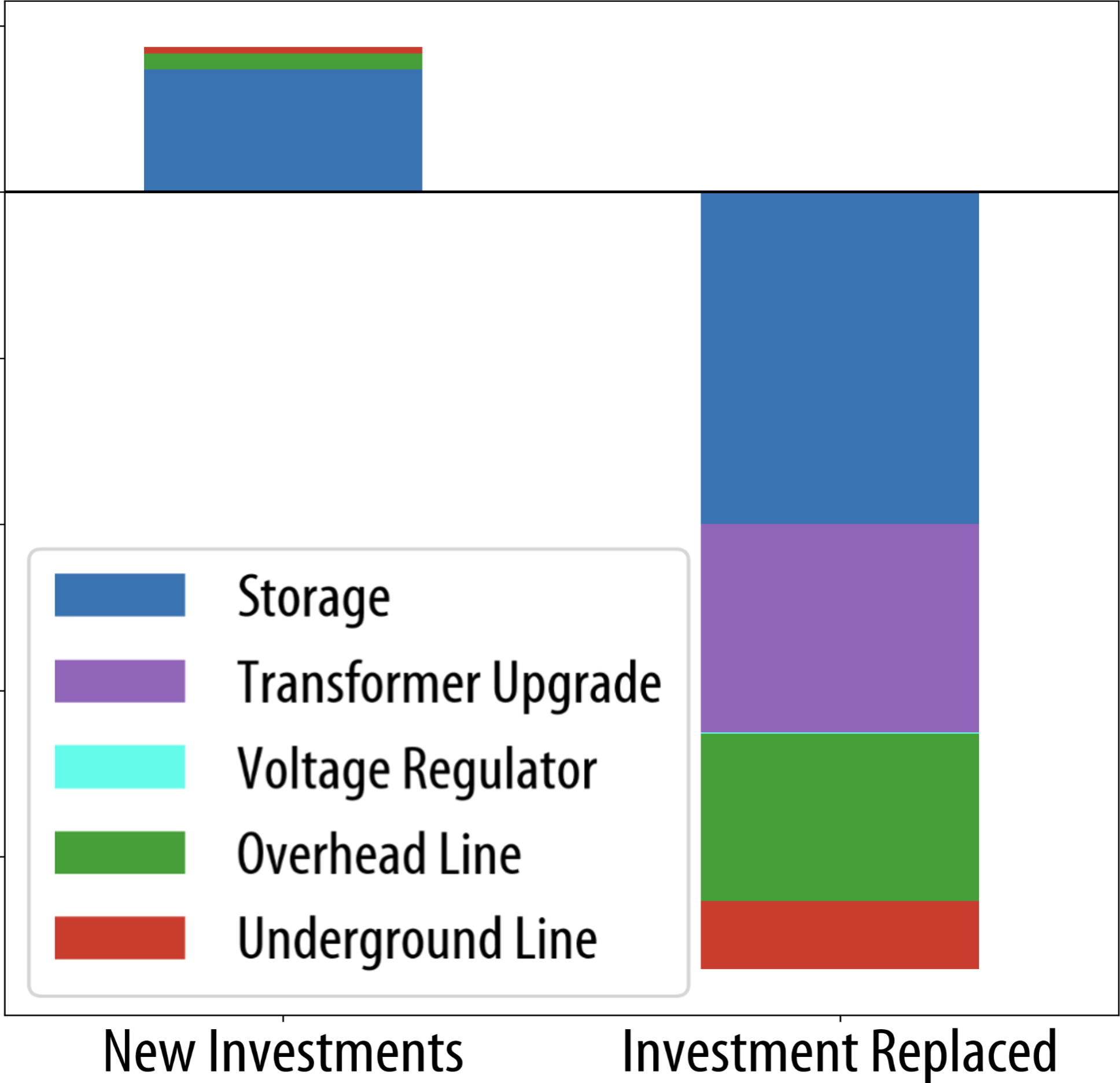}
        \label{fig:incremental_deferrals}
    }
        \subfigure[Deferral type per feeder]{
        \centering
        \includegraphics[width=0.22 \textwidth]{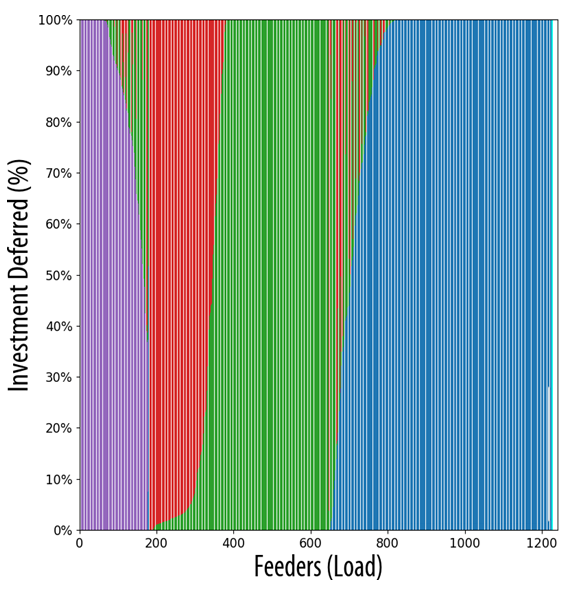}
        \label{Fig:deferral_feeders}
    }
      \caption{Incremental costs ("High Rooftop PV" scenario) and deferrals ("High Load" scenario) per type of upgrade and feeder. “New Investments” and “Investments Replaced” refer to the upgrades with and without CS,respectively.}
  \label{fig:invesmtent_type}
\end{figure}

\begin{figure}[!ht]
       \centering
       \includegraphics[width=0.45 \textwidth]{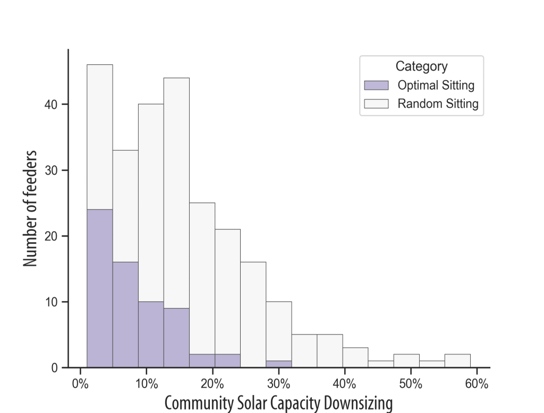}
      \caption{The number of feeders with CS capacity downsizing as an economically viable solution for CS integration. Comparison between random and optimal siting (PV Case)}
      \label{Fig:siting and downsizing}
\end{figure}

\begin{figure*}[!ht]
\centering
\includegraphics[width=0.95\linewidth]{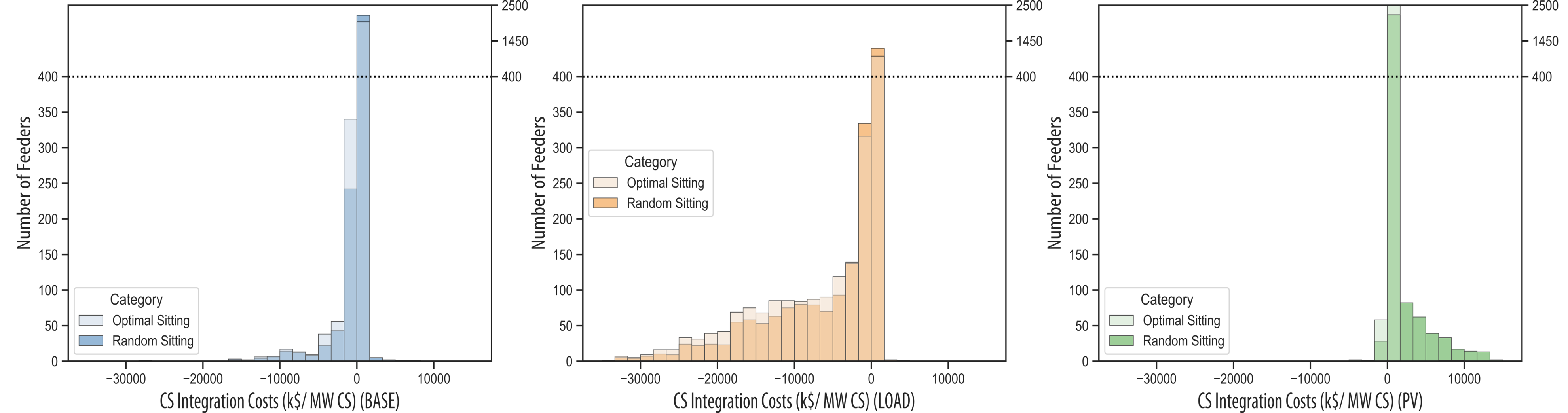}
\caption{Impact of strategic siting on overall CS integration costs in the three scenarios.}
\label{fig:strategic_st}
\end{figure*}

When normalizing results relative to the CS size, as depicted in Figures~\ref{fig:cost_cs} and \ref{Fig: deferrals_cs}, the grid costs and benefits associated with the presence of CS are generally below \$50/kW installed for the majority of feeders. Assuming CS costs fall within the range of \$1.5k/kW, these integration costs amount to approximately 3\% of this value, which may not pose significant limitations to development. However, for feeders approaching their hosting capacity limit, costs can exceed \$250/kW of CS installed in specific cases. Conversely, in heavily loaded feeders, the benefits may surpass \$500/kW of CS installed in certain instances.

When quantifying these values from the perspective of consumers in Figures~\ref{fig:cost_cons} and \ref{Fig: deferrals_cons}, CS-related costs and deferral benefits are generally moderate for most feeders, around 0.2~\cent/kWh. However, in exceptional circumstances, these costs and benefits could escalate to 2\cent~per kWh of feeder energy demand. Interestingly, this 2~\cent/kWh figure also serves as a benchmark for the component of total distribution grid costs in consumer electricity rates \cite{osti_1840705}. Consequently, in extreme scenarios, CS integration could represent 100\% of the total distribution costs. Nevertheless, under normal circumstances, this value tends to be around 10\%.

\subsection{Investment types}

Figure \ref{fig:invesmtent_type} illustrates a breakdown of investments and deferrals categorized by the type of upgrades. In cases where CS integration requires additional costs, these are predominantly associated with new investments in overhead and underground reconductoring, as depicted in Figures~\ref{fig:incremental_investments} and \ref{Fig:investments_feeder}). 

The incremental deferrals in Figure~\ref{fig:incremental_deferrals} show that, even when  are integration costs ($C^{itgr}_{CS}$) are negative, upgrades may still be be necessary. However, these investments are observed to replace a significantly larger amount of upgrades. Another relevant aspect, evident in both Figures ~\ref{fig:incremental_deferrals} and \ref{Fig:deferral_feeders}, is the diversity of deferred upgrade types, encompassing a mix of transformer, storage, and reconductoring.

Lastly, in Figure ~\ref{fig:incremental_deferrals}, it is interesting to observe that utility-owned storage can function both as a deferral and a deferred investment, depending on the specific case.

\subsection{Strategic siting and capacity downsizing}

To analyze the impact of strategically siting CS, we ran the cost assessment for two cases: a) when CS is randomly placed within the pre-defined locations, top, middle, and bottom of the feeder. b) when CS location is a decision variables of the optimal capacity expansion model, as discussed above.

When grid infrastructure upgrades are neither technically feasible nor economically viable for CS integration, the model may choose to curtail CS generation. This decision has a cost to the system (assumed equal to the locational marginal price) and, in practice, it implies a downsizing of the CS capacity. 

Figure \ref{Fig:siting and downsizing} shows that, even when CS is randomly located, capacity downsizing is infrequent, occurring only in 359 cases, most of them with capacity downsizing needs below 20\%. When strategic siting is considered, and CS is optimally sited to reduce integration costs, the downsizing needs to reduce even further and become marginal (only in 66 feeders).

Furthermore, Figure \ref{fig:strategic_st} shows that strategic siting of CS is also effective in reducing upgrade costs and enhancing deferral value across the three scenarios studied. It is interesting to note that the distribution of costs is shifted to the left when optimal siting is considering. These results suggest that CS projects facing high integration costs should explore alternate sites on the same feeder as a cost reduction strategy.

\subsection{Policy and regulatory takeaways}

Under current feeder conditions, most CS projects do not require additional infrastructure upgrades: only 7\% of the feeders required upgrades to support CS in the base case. This shows that the hosting capacity benchmark for MDL to size CS may be too conservative. Regulators could require utilities to perform more rigorous CS interconnection studies to identify larger sizes for CS.
In contrast, downsizing CS projects due to distribution system constraints is rarely necessary, but we find that small capacity downsizing (below 20\%) may be cost-effective for 10\% of feeders Regulatory. CS requirements could include CS project downsizing guidelines to reduce interconnection costs for CS developers.

When a CS project results in net deferral benefits, it can still require some distribution system upgrades. This suggests that interconnection analyses that focus exclusively on grid upgrades and not on deferrals may overestimate CS interconnection costs. Regulatory frameworks should encompass quantification of both costs and potential deferrals, with and without CS, to fully assess the economic impacts of CS integration. Another avenue for capturing deferral benefits in regulation is to require utilities to identify feeders that are candidates for non-wires alternatives, so CS developers can target these feeders for project siting. Shared-savings mechanisms could incentivize utilities to work with CS developers to achieve capacity deferral benefits on those feeders, while providing savings to utility customers. 

Strategic siting within a feeder can substantially reduce interconnection costs and reduce the need for downsizing CS projects. This suggests that joint evaluations between CS developers and utilities to find optimal locations for project siting could be beneficial. 
Utilities could be required to work with developers to identify feeder locations that are more strategic for siting CS projects.

\section{Conclusions}\label{sec.Conclusions}

This work provided a comprehensive study on distribution grid capital investments required to integrate CS at the feeder level. The methodology, designed to evaluate different scenarios and compare feeders with different sizes, relied on incremental least-cost expansion of the distribution system, which allowed to capture not only feeder investments but also deferrals and to identify specific trends on equipment upgrade types related to CS integration. This methodology was applied to more than 2,500 representative feeders, with 3 loading conditions and the option of strategic siting, in a total of more than 15,000 instances, and compared with empirical costs from 210 real CS integration projects.

We found that current practices of CS sizing result in minor or no distribution upgrades at all in most cases. When upgrades are required, reconductoring is often the most significant investment, and mitigation options such as strategic siting and downsizing of CS can help reduce them. Additionally, the potential for deferrals is at least as high as investment costs, which indicates that specific regulatory possesses for CS integration should start factoring in those benefits.

Future works could leverage the costs and benefits quantified in this work to study potential policy and regulatory mechanisms for cost allocation of CS related upgrades and discuss economic impacts in CS subscribers and non-subscribers.

\section*{Acknowledgment}
This work was funded by the U.S. Department of Energy under the contract DE-AC02-05CH11231. We also thank our colleagues Greg Leventis, Chandler Miller and Jeff Deason for the valuable discussion during this work.

\bibliographystyle{IEEEtran}
\bibliography{IEEEabrv,References}

\end{document}